\documentclass[12pt,leqno]{amsart} 
\usepackage[numbers]{natbib}
\usepackage{amsmath,amsthm}
\usepackage{amssymb}
\usepackage{stmaryrd}
\renewcommand{\widehat}{\hat}
\usepackage{url}
\usepackage[letterpaper,colorlinks={true},
backref,pdfstartview={FitBH},pdfview={FitBH},bookmarks={true}]{hyperref}
\global\hbadness=500
\global\tolerance=200
\newtheorem{theorem}{Theorem}[section]
\newtheorem{lemma}[theorem]{Lemma}
\newtheorem{conjecture}[theorem]{Conjecture}

\newtheorem{question}[theorem]{Question}

\newtheorem{corollary}[theorem]{Corollary} 
\theoremstyle{definition}
\newtheorem{definition}[theorem]{Definition}

\theoremstyle{remark}

\numberwithin{equation}{theorem}

\newcommand{\ce}{computably enumerable }
 
\newcommand{\join}{\oplus} 
\renewcommand{\phi}{\varphi}

\newcommand{\E}{\mathcal{E}}
\renewcommand{\L}{\mathcal{L}}\renewcommand{\S}{\mathcal{S}}

\newcommand{\F}{\mathcal{F}}
\newcommand{\R}{\mathcal{R}}
\newcommand{\B}{\mathcal{B}}
\newcommand{\D}{\mathcal{D}}
\newcommand{\C}{\mathcal{C}}

\newcommand{\Ahat}{\widehat{A}} \newcommand{\Bhat}{\widehat{B}}

\newcommand{\Dhat}{\widehat{D}}

\newcommand{\Rhat}{\widehat{R}}
\newcommand{\Hhat}{\widehat{H}}
\newcommand{\What}{\widehat{W}}

\newcommand{\That}{\widehat{T}}
\newcommand{\Shat}{\widehat{S}}

  \newcommand
{\wock}{\omega_1^{\textup{CK}}}

\title[Orbits]{The Complexity of  Orbits of Computably\\ Enumerable Sets}
  
\author[P.\ Cholak]{Peter~A.~Cholak}

\address{Department of Mathematics\\ University of Notre Dame\\ 
  Notre Dame, IN 46556-5683}

\email{Peter.Cholak.1@nd.edu}

\urladdr{http://www.nd.edu/\~{}cholak}

\author[R.\ Downey]{Rodney Downey}

\address{School of Mathematics, Statistics and Computer Science\\
 Victoria University\\ P.O. 
Box 600\\ Wellington, New Zealand}

\email{Rod.Downey@vuw.ac.nz}
\urladdr{http://www.mcs.vuw.ac.nz/\~{}downey} 

\author[L.\ Harrington]{Leo~A.~Harrington}

\address{Department of Mathematics\\ University of California \\
  Berkeley, CA 94720-3840}

\email{leo@math.berkeley.edu}

\date{\today}

\thanks{Research partially supported NSF Grants DMS-96-34565,
  99-88716, 02-45167 (Cholak), Marsden Fund of New Zealand (Downey),
  DMS-96-22290 and DMS-99-71137 (Harrington).  Some of involved work
  was done partially while Cholak and Downey were visiting the
  Institute for Mathematical Sciences, National University of
  Singapore in 2005. These visits were supported by the Institute.}

\begin{document}


\begin{abstract}
  The goal of this paper is to announce there is a single orbit of the
  c.e.\ sets with inclusion, $\E$, such that the question of
  membership in this orbit is $\Sigma^1_1$-complete.  This result and
  proof have a number of nice corollaries: the Scott rank of $\E$ is
  $\wock +1$; not all orbits are elementarily definable; there is no
  arithmetic description of all orbits of $\E$; for all finite $\alpha
  \geq 9$, there is a properly $\Delta^0_\alpha$ orbit (from the
  proof).
\end{abstract}

\maketitle

\section{Introduction}

In the classic paper \cite{Post:44}, Post suggested that the study of
the lattice of computably (recursively) enumerable (c.e.) sets was
fundamental in computability theory.  Post observed that, at the time,
all known undecidability proofs worked by coding some
``noncomputability'', (coded by a certain kind of set) into the theory
at hand thereby arguing that the relevant structures could emulate
computation.  He argued that for such proofs the central object was
the notion of a computably enumerable set, being one which is
effectively generated as the range of some effective process.  The
basic example is the set of consequences of a computably enumerable
set of axioms for a formal system.  Of course, the other key concept
in computability was that of Turing \cite{Turing:39} who introduced
the notion of reducibility.  Reducibilities are pre-orderings used to
measure relative computational complexity.

The interplay of these two basic objects, (Turing) reducibility and
effectively enumerated (c.e.) sets has a long and rich history.
Clearly the computably enumerable sets under union and intersection
form a lattice, denoted by ${\mathcal E}$.  Their Turing degrees form
an upper semilattice, denoted by ${\mathcal R}$.  Ever since the
ground breaking paper of Post, there has been a persistent intuition
that structural properties of computably enumerable sets have
reflections in their degrees, and conversely.  In particular, {\em
  definability} in ${\mathcal E}$ should be linked with {\em
  information content} as measured by ${\mathcal R}.$

The simplest possible illustration of this is the fact that the
complemented members of ${\mathcal E}$ are exactly the members of
${\bf 0}$ the degree of the computable sets.  An excellent and deep
example is Martin's result that the Turing degrees of maximal sets are
exactly the high \ce Turing degrees\footnote{Indeed the reader should
  recall that, more generally, a set $A$ is low$_n$ iff $A^{(n)} =
  \emptyset^{(n)}$ iff $\Delta^0_{n+1} = \Delta^A_{n+1}$, and $A$ is
  high$_n$ iff $A^{(n)} = \emptyset^{(n+1)}$ iff $\Delta^0_{n+2} =
  \Delta^A_{n+1}$.}, (that is, their jumps are as complex as possible,
${\bf a'}={\bf 0''}$), where a co-infinite set $A$ is a maximal c.e.\
set iff for all c.e.\ sets $B$, if $A\subseteq B$ then either $A=^*B$
or $B=^* \omega$, where $=^*$ denotes equivalence modulo the filter of
finite sets.  Since a set $A$ is finite iff every subset is
complements in ${\mathcal E}$, it is natural to consider
$\mathcal{E}^*$, $\mathcal{E}$ modulo the filter of finite sets.  Thus
$A$ is a maximal set iff it represents a co-atom in ${\mathcal E^*}.
$

An original impetus for work on ${\mathcal E^*} $ was {\em Post's
  Problem} : Post observed that the coding inherent in all known
undecidability proofs of the time were so faithful that all computably
enumerable sets seemed to be either the ${\bf 0'}$ the degree of the
halting problem, or ${\bf 0}$ the degree of the computable sets. He
asked if this was always the case. Post's problem is the following :
{\em Are there c.e.\ sets of intermediate Turing degree? That is there
  a \ce degree ${\bf a }$ with ${\bf 0}<{\bf a}<{\bf 0'}.$} Post even
suggested a programme to answer this question.  Since complemented
members of ${\mathcal E}$ are computable, but
$\emptyset'=\{x:\varphi_x(x)\downarrow\}$, the halting problem, had
many infinite computably enumerable sets disjoint from it, perhaps a
very ``noncomplemented'' member of ${\mathcal E}$ would turn out to be
incomplete. Towards this goal, Post defined a computably enumerable
set $A$ to be {\it simple} if $\omega-A$ was infinite and for all
infinite computably enumerable $W$, $W\cap A\neq \emptyset$. Post
constructed a simple set and then proved that a simple set $A$ is not
of the same $m$-degree as ${\bf 0'}$.  Post's programme was to come up
with a thinness property of the complement of a c.e.\ set which would
guarantee Turing incompleteness.  Interpreted more liberally, we can
generalize Post's programme and ask whether there is any definable
property of a c.e.\ set in ${\mathcal E^*}$ which would guarantee
Turing incompleteness.

While Post's problem was eventually solved by the development of the
priority method independently by Friedberg \cite{Friedberg:57} and
Muchnik \cite{Muchnik:56} out of the work of Kleene and Post
\cite{Kleene.Post:54}, whether Post's Programme could be realized
successfully remained open for much longer.  Myhill observed that if
Post's original programme were to succeed then {\it maximal sets},
first constructed by Friedberg, should be Turing incomplete.  On the
other hand, rather than being Turing incomplete, Tennenbaum suggested
that all maximal sets would be Turing {\em complete}.  As we have have
seen above, following earlier work of Tennenbaum, Sacks, and Yates,
and others, Martin classified the degrees of maximal sets as precisely
the high \ce degrees.  Martin's theorems demonstrate can be seen as a
partial realization of Tennenbaum's intuition.  That is, in spite of
the fact that maximal sets may not necessarily be complete, they do
have high information content in the sense that they have the same
jump as the halting problem.  That is, as far as the jump operator is
concerned, they are indistinguishable from the halting problem.  Martin
had discovered the first {\it invariant class} in ${\mathcal R}$ in
the sense that the high degrees are precisely those realized by the
(definable class of) maximal sets.

Definability goes hand in hand with automorphisms of structures.
Thus, a class ${\mathcal C}$ of computably enumerable degrees {\it
  invariant} if there is a collection ${\mathcal C}'$ of computably
enumerable sets closed under automorphisms of ${\mathcal E}$ such that
${\mathcal C}=\{\mbox{deg}(A):A\in {\mathcal C}'\}$.  The following
definition will be important for our subsequent discussion.

\begin{definition}
  $A \approx \Ahat$ iff there is a map, $\Phi$, from the c.e.\ sets to
  the c.e.\ sets preserving inclusion, $\subseteq$, (so $\Phi \in
  \text{Aut}(\E)$) such that $\Phi(A)=\Ahat$.
\end{definition}

By \citet{Soare:74}, $\mathcal{E}$ can be replaced with
$\mathcal{E}^*$, since Soare showed that every automorphism of
${\mathcal E^*}$ is equivalent to one on ${\mathcal E}$ and conversely
(as long we focus on infinite and co-infinite sets).  Early work by
Lachlan and others showed that both of the automorphism groups
Aut(${\mathcal E})$ and Aut$({\mathcal E}^*)$ were large since each
had $2^{\aleph_0}$ automorphisms. Martin \cite{Martin:66*1} used a
priority construction to show that a certain construction of Post
(hypersimplicity) was not invariant under automorphisms of ${\mathcal
  E}$.

Post's original programme was to look at thinness properties of the
lattice of the complement of a c.e.\ set $A$. We will see shortly this
original programme cannot be solved. But there several solutions to
various modified versions of this programme.

The earliest solution to a modified Post's Programme was due to
Marchenkov \cite{Marchenkov:76} who showed that a certain type of
maximal set in a related quotient structure gave a solution.
Specifically, if you change the game and replace the integers by
computably enumerable equivalence classes $\eta$, you can get a
structure consisting of the c.e.\ sets factored out by this
equivalence relation. Then for a suitable choice of $\eta$ certain
$\eta$-maximal sets provide a solution to Post's programme, at least
in a generalized
sense. 

However we want to focus on the structures $\mathcal{E}$ and
$\mathcal{E}^*$. In these structures, there is a positive solution to
another modified version of Post's Programme.  This solution will be
discussed shortly but for now we want to focus on the failure of the
original programme. This leads us to a definition which will be
important:

\begin{definition}\label{la}
  $\L^*(A)$ is $ \{ W \cup A : W \text{ an c.e.\ set}\}$ under
  $\subseteq$ modulo the ideal of finite sets ($\F$). (The outside of
  a set.)
\end{definition}

The final blow to Post's original programme was the marvelous paper of
Soare \cite{Soare:74}, who showed that maximal sets form an orbit in
Aut$({\mathcal E})$.  In particular, no ``extra'' property together
with maximality could guarantee incompleteness.  Moreover, the paper
Cholak, Downey and Stob \cite{Cholak.Downey.ea:92*1}, showed that no
property of $\L^*(A)$ {\em alone} could guarantee Turing
incompleteness for a \ce set $A$.  That is, Cholak, Downey and Stob
proved that if for any \ce set $A$ there is a set \ce Turing complete
$B$ with the same lattice of supersets as $A$, $\L^*(A)\cong \L^*(B)$.

Soare's \cite{Soare:74} was highly influential. The methods introduced
constructed $\Delta_3^0$ automorphisms of ${\mathcal E^*}$.  Here we
will classify automorphisms according to the complexity of a
presentation of them. A presentation will be a function
$f:\omega\to\omega$ such that that $W_e\mapsto W_{f(e)}$ induces the
automorphism.  If $e\mapsto f(e)$ a $\Delta_3^0$ function, we would
call the automorphism $\Delta_3^0$, and $A\approx_{\Delta_3^0}
\Ahat$. While later papers presented Soare's automorphism machinery
argument as a more thematic and flexible tree argument (beginning with
Cholak \cite{Cholak:91} and \cite{Cholak:95} and Harrington and Soare
\cite{Harrington.Soare:96}) most of the key underlying ideas for
constructing automorphisms of $({\mathcal E})$ are in Soare's original
paper.

The principal tool used is called the (or, in view of recent work,
{\em an}) {\em Extension Lemma}.  Roughly speaking, constructing an
automorphism works as follows.  We wish to show maximal $A \approx
\Ahat$.  We are given two versions of the universe $\omega$, called
$\omega$ and $\hat{\omega}$ with $A\subset \omega$ and $\Ahat \subset
\hat{\omega}$, with enumerations of c.e.\ $\{W_e:e\in \omega\}$
subsets of $\omega$ and $\{V_e:e\in \omega\}$ subsets of
$\hat{\omega}$.  We must define some sort of mapping as follows

\begin{center}
  $
  \begin{array}{lll}
    \omega & \mapsto &\hat{\omega}\\
    A&\mapsto&\hat{A}
    \\
    W_e&\mapsto
    &
    \hat{W}_{f(e)} \\
    \hat{V}_{g(e)}&\mapsfrom &V_{e}.\\
  \end{array}$

\end{center}

We must have the $\hat{W}_{f(e)}$ and $\hat{V}_{g(e)}$ so that we can
argue that the mapping induces an automorphism of $\E^*$ by a back and
forth argument.  Here we are thinking of {\em building} the hatted
sets $\hat{W}$ and $\hat{V}$.  At the very least, intersections should
be respected.  That is, if $\overline{A}\cap W_e$ is infinite, then
$\overline{\Ahat}\cap \hat{W}_{f(e)}$ would need to be infinite.
Similarly if $W_e-\hat{V}_{g(q)}$ is infinite the so too must be
$\hat{W}_{f(e)}-V_{q}, $ etc.  Evidently, any possible diagram we can
think of denoting intersections and difference would need to be
respected.  We represent these intersections and differences by {\em
  states}, which are strings measuring which $W_e$'s and $\hat{V}_k$'s
on the $\omega$ side an element is in, and which $\hat{W}_j$'s and
$V_q$'s a hatted element is in on the $\hat{\omega}$ side.  We would
write this basic requirement as
$$\R_\sigma:\exists ^\infty x \in \mbox{ state }\sigma\mbox{ iff }
\exists^\infty \hat{x} \in \mbox{ state }\hat{\sigma}.$$

Soare's original idea is to begin on the $\L^*(A)\cong \L^*(\hat{A}).$
He would make this an isomorphism and then extend this outside
isomorphism to an automorphism by an isomorphism of the lattice of \ce
subsets of $A$ to those of $\hat{A}$.

Concentrating on the $\L^*(A)\cong \L^*(\hat{A})$ part, as we go
along, elements appear to be in this region (that is, in
$\overline{A}_s$ or $\overline{\hat{A}}_s$), and we build
corresponding sets to match the states measuring intersections.  This
would seem not too hard in the case of a maximal set since for any
sets $W$ either $W\cap \overline{A}$ is finite or $W$ almost contains
$\overline{A}$. On the hatted side, all we would need to do is either
have the corresponding $\hat{W}$ empty, or containing
$\overline{\Ahat}$, and similarly for the mappings from the hatted
side back.  The information as to which is correct is $\Sigma_3^0$
information and can be handled by a priority argument.

However, the {\em heart} of Soare's method is the following.  As we go
along enumerating hatted sets as elements stream into $\overline{A}_s$
many of these will be based on {\em wrong} information (such as the
fact that at stage $s$ they might appear in the complement of $A$ yet
might be in $A$), and will later enter $A$, and hence be in
$W_e\searrow A$.  They will enter $A_t$ and $\Ahat_u$ in various {\em
  entry} states some caused by the $W_e$ and $V_e$ played by the
opponent, and the hatted sets played by {\em us}.  The key problem is
how to handle these wrongly enumerated elements and be able to {\em
  extend} the $\L^*(A)\cong \L^*(\hat{A})$ correspondence to an
automorphism.  For the automorphism machinery to succeed, it is {\it
  necessary} to ensure that {\it for all {\em entry} states $\sigma$
  (and dually for $\widehat{\sigma}$), if infinitely many elements
  enter $A$ in state $\sigma$, then there is some covering entry state
  $\widehat{\tau}$ ($\tau$, respectively).}  Here {\em covering} means
that it is within {\em our} power to add elements into sets under {\em
  our} control to be able to match states.

Soare's Extension Lemma shows that this necessary condition is
sufficient.  Soare showed that we have not already killed that
automorphism, meaning that the necessary condition is satisfied, then
there is a strategy which enables us to extend the partial matching
into a full automorphism.  On the inside, that is the lattice of
subsets of $A$ to those of $\Ahat$ the map is $\Delta_3^0$, and in the
case of maximal sets, Soare's original result shows that if $A$ and
$\Ahat$ are maximal, then $A\approx_{\Delta_3^0} \Ahat.$

There is a lot of subsequent work on automorphisms and invariance in
the lattice of c.e.\ sets.  Almost all of it either uses Soare's
original Extension Lemma as a black box, or modified it, to prove
various results on the lattice of c.e.\ sets.  Examples include the
work of Maass \cite{Maass:84}, Maass and Stob \cite{Maass.Stob:83},
and Downey and Stob \cite{Downey.Stob:92}.

Early on, the methods seemed so powerful that anything seemed
possible. Perhaps all sets were automorphic to complete sets, as
suggested by Soare \cite{Soare:74*1}. Certainly Harrington and Soare,
and Cholak independently showed that all sets were automorphic to high
sets.

Hand in hand with this work constructing automorphisms was another
line of investigation, where {\em failures} of the automorphism
machinery could be exploited to provide {\em definability} results in
$\mathcal{E}^*$.  A classic example of this is the following theorem
of Harrington and Soare \cite{Harrington.Soare:91}
\cite{Harrington.Soare:96*1} who showed that a more general form of
Post's Programme indeed has a positive solution.


\begin{theorem}[Harrington and Soare \cite{Harrington.Soare:91}]
  There is a definable property $Q(A)$, such that, if a c.e.\ set $A$
  satisfies $Q(A)$, then $A$ is Turing incomplete.
\end{theorem}

There were precursers to the Harrington-Soare result.  Harrington used
the idea of exploiting the failure of the machinery to get a
definition of being a halting problem in the lattice of c.e.\
sets. Similarly Lerman and Soare \cite{Lerman.Soare:80*1} showed that
there are low simple sets that are elementarily inequivalent, in that
one has a property called d-simplicity and one has not, where
d-simplicity is an elementary property implying certain facts about
entry states.  Another example of this can be found in Downey and
Harrington \cite{Downey.Harrington:96} where the ``no fat orbit''
theorem is proven.  The simplest form of the Downey-Harrington result
below says that no c.e.\ set has an orbit hitting all nonzero degrees.

\begin{theorem}[Downey and Harrington -- No fat orbit]
  There is a property $S(A)$, a prompt low degree $\mathbf{d_1}$, a
  prompt high$_2$ degree $\mathbf{d_2}$ greater than $\mathbf{d_1}$,
  and tardy high$_2$ degree $\mathbf{e}$ such that for all $E \leq_T
  \mathbf{e}$, $\neg S(E)$ and if $\mathbf{d_1} \leq_T D \leq_T
  \mathbf{d_2}$ then $S(D)$.
\end{theorem}

We remark that this ``failure'' methodology has yielded similar
definability results in other structures such as the lattice of
$\Pi_1^0$ classes, as witnessed by Weber \cite{Weber:04} and
\cite{Weber:06}, Cholak and Downey \cite{MR2005e:03092}, and Downey
and Montalb{\'a}n \cite{Downey.Montalban:nd}.  Perhaps the best
example of the methodology is the following proof of the definability
of the double jump classes, the proof using ``patterns'' which are
more or less direct reflections of blockages to the automorphism
machinery.

\begin{theorem} [Cholak and Harrington 02]
  Let $\mathcal{C} = \{\mathbf{a}: \mathbf{a}$ is the Turing degree of
  a $\Sigma_3 \text{ set greater than } \mathbf{0''}\}$.  Let
  $\mathcal{D} \subseteq \mathcal{C}$ such that $\mathcal{D}$ is
  upward closed.  Then there is an non-elementary
  ($\L_{\omega_1,\omega}$) $\mathcal{L}(A)$ property
  $\varphi_{\mathcal{D}}(A)$ such that $D''\in \mathcal{D}$ iff there
  is an $A$ where $A \equiv D$ and $\varphi_{\mathcal{D}}(A)$.
\end{theorem}

\begin{corollary}
  If $\mathbf{a}'' > \mathbf{b}''$ then there is a $A \in \mathbf{a}$
  such that for all $B\in \mathbf{b}$, $A$ is not automorphic to $B$
  (in fact, $\mathcal{L}^*(A) \not\cong \mathcal{L}^*(B)$).
\end{corollary}

Related here is the following conjecture of Harrington.

    \begin{conjecture}[Harrington]
      For all $A$ and degrees $\mathbf{d}$ if $A' \leq_T \mathbf{d}'$
      is there $\Ahat \in \mathbf{d}$ such that $\mathcal{L}^*(A)
      \cong \mathcal{L}^*(\Ahat)$.
    \end{conjecture}

    For more of these results one can see the paper
    \cite{Cholak.Harrington:00}.

    \section{New Results}

    The present work is motivated by basic questions about the
    automorphism group of ${\mathcal E}^*$.  How complicated is it? If
    $A\approx\Ahat$ is $A\approx \Ahat$ witnessed by an {\em
      arithmetical} automorphism?  How complicated is
    $\{W_e:W_e\approx A\}$ for a fixed $A$?  The following conjecture
    was made by Ted Slaman and Hugh Woodin in 1989.

\begin{conjecture}[\citet{Slaman.Woodin:conjecture}]
  The set $\{\langle i , j \rangle: W_i\approx W_j)\}$ is
  $\Sigma^1_1$-complete.
\end{conjecture}

This conjecture was claimed to be true by the authors in the mid
1990s; but no proof appeared.  One of the roles of this announcement
and the full paper \cite{cholak.downey.ea:nd} is to correct that
omission. The proof we will present is far simpler than all previous
(and hence unpublishable) proofs.  Indeed, much of the material
reported in this paper due to Cholak and Harrington was developed
towards making the proof of the conjecture accessible. The other
important role of this communcation and the full paper is to prove a
stronger result.

\begin{theorem}[The Main Theorem]\label{main}
  There is a c.e.\ set $A$ such that the index set $\{i : W_i \approx
  A\}$ is $\Sigma^1_1$-complete.
\end{theorem}
   
As mentioned in the abstract this theorem does have a number of nice
corollaries.

\begin{corollary}
  Not all orbits are elementarily definable; there is no arithmetic
  description of all orbits of $\E$.
\end{corollary}

\begin{corollary}
  The Scott rank of $\E$ is $\wock +1$.
\end{corollary}

\begin{proof}
  Our definition that a structure has Scott rank $\wock +1$ is that
  there is an orbit such that membership in that orbit is
  $\Sigma^1_1$-complete. There are other equivalent definitions of a
  structure having Scott Rank $\wock +1$ and we refer the readers to
  \citet{MR1767842}.
\end{proof}

A consequence of the {\em method} of the proof (and some further
effort to preserve quantifiers) is the following.

\begin{theorem}\label{sec:maincor}
  For all finite $\alpha > 8$ there is a properly $\Delta^0_{\alpha}$
  orbit.
\end{theorem}

Hitherto this paper \cite{cholak.downey.ea:nd} all known orbits were
$\Delta_3^0$ with the single exception of the orbit of Cholak and
Harrington \cite{Cholak.Harrington:nd} which constructed a pair of
sets $\Delta_5^0$ automorphic but not $\Delta_3^0.$
  
Before we turn to the proof of Theorem~\ref{main}, we will discuss the
background to the Slaman-Woodin Conjecture.  Certainly the set
$\{\langle i , j \rangle: W_i\approx W_j)\}$ is $\Sigma^1_1$.  Why
would we believe it to be $\Sigma_1^1$-complete?  The following result
is from the folklore\footnote{We think it is well known that the
  isomorphism problem for Boolean Algebras and Trees are
  $\Sigma^1_1$-complete, at least in the form stated in
  Theorems~\ref{folkloreBA} and \ref{folkloreI}. We have searched for
  a reference to a proof for these theorems without success.  It seems
  very likely that these theorems were known to Kleene. There are a
  number of places where something very close to what we want appears;
  for example, see the example at the end of Section~5 of
  \citet{MR2058190} and surely there are earlier examples (for
  example, \citet{walkerwhite:00}).  All of these constructions work
  by coding the Harrison ordering. In the full paper we give
  self-contained proofs of the folklore theorems we use.}.

\begin{theorem}[Folklore\footnote{See Section~5 of the full paper
    \cite{cholak.downey.ea:nd} for more information and a
    proof.}] \label{folkloreBA} There is a computable listing,
  $\mathcal{B}_i$, of computable Boolean algebras such that the set
  $\{\langle i , j \rangle: \mathcal{B}_i\cong \mathcal{B}_j\}$ is
  $\Sigma^1_1$-complete.
\end{theorem}

Note that $\L^*(A)$ is a definable structure in $\E$ with a parameter
for $A$.  The following result says that the full complexity of the
isomorphism problem for Boolean algebras of Theorem \ref{folkloreBA}
is present in the supersets of a c.e.\ set.

\begin{theorem}[\citet{Lachlan:68*3}]
  Effectively in $i$ there is a c.e.\ set $H_i$ such that $\L^*(H_i)
  \cong \mathcal{B}_i$.
\end{theorem}

\begin{corollary}
  The set $\{\langle i , j \rangle: \L^*(H_i)\cong \L^*(H_j)\}$ is
  $\Sigma^1_1$-complete.
\end{corollary}
  
Slaman and Woodin's idea was to replace ``$\L^*(H_i)\cong \L^*(H_j)$''
with ``$H_i\approx H_j$''.  Unfortunately, this very attractive idea
is doomed, as we now see.

\begin{definition}[The sets disjoint from $A$]
  $$\mathcal{D}(A) = ( \{ B : \exists W ( B \subseteq A \cup W \text{ and
  } W \cap A =^* \emptyset)\}, \subseteq).$$ Let
  $\mathcal{E}_{\mathcal{D}(A)}$ be $\mathcal{E}$ modulo
  $\mathcal{D}(A)$.  $A$ is \emph{$\mathcal{D}$-hhsimple} iff
  $\mathcal{E}_{\mathcal{D}(A)}$ is a Boolean Algebra.  $A$ is
  \emph{$\mathcal{D}$-maximal} iff $\mathcal{E}_{\mathcal{D}(A)}$ is
  the trivial Boolean Algebra.
\end{definition}

\begin{lemma}
  If $A$ is simple then $\mathcal{E}_{\mathcal{D}(A)}
  \cong_{\Delta^0_3} \L^*(A)$.
\end{lemma}

It is an old result of Lachlan \cite{Lachlan:68*3} that $A$ is
hhsimple iff $\mathcal{E}_{\mathcal{D}(A)}$ is a Boolean
algebra. Except for the creative sets, until recently, all known
orbits were orbits of $\mathcal{D}$-hhsimple sets.  We direct the
reader to \citet{Cholak.Harrington:nd} for a further discussion of
this claim and for an orbit of $\E$ which does not contain any
$\mathcal{D}$-hhsimple sets.  The following are relevant theorems from
\citet{Cholak.Harrington:nd}.

\begin{theorem}\label{dhhsimple}
  If $A$ is $\mathcal{D}$-hhsimple and $A$ and $\Ahat$ are in the same
  orbit then $\E_{\mathcal{D}(A)} \cong_{\Delta^0_3}
  \E_{\mathcal{D}(\Ahat)}$.
\end{theorem}

\begin{theorem}[using \citet{Maass:84}]
  If $A$ is $\mathcal{D}$-hhsimple and simple (i.e., hhsimple) then
  $A\! \approx \Ahat$ iff $\L^*(A) \cong_{\Delta^0_3} \L^*(\Ahat)$.
\end{theorem}

Hence the Slaman-Woodin plan of attack fails.  In fact even more is
true.

\begin{theorem}
  If $A$ and $\Ahat$ are automorphic then
  $\mathcal{E}_{\mathcal{D}(A)}$ and
  $\mathcal{E}_{\mathcal{D}(\Ahat)}$ are $\Delta^0_6$-isomorphic.
\end{theorem}

Hence in order to prove Theorem~\ref{main} we must code everything
into $\D(A)$.  This is completely contrary to all approaches used to
try to prove the Slaman-Woodin Conjecture over the years.  We will
point out two more theorems from \citet{Cholak.Harrington:nd} to show
how far the sets we use for the proof must be from simple sets, in
order to prove Theorem~\ref{main}.

\begin{theorem}
  If $A$ is simple then $A \approx \Ahat$ iff $A \approx_{\Delta^0_6}
  \Ahat$.
\end{theorem}

\begin{theorem}
  If $A$ and $\Ahat$ are both promptly simple then $A \approx \Ahat$
  iff $A \approx_{\Delta^0_3} \Ahat$.
\end{theorem}

\section{Future Work and the Degrees of the Constructed Orbits}

While this work does answer many open questions about the orbits of
c.e.\ sets, there are many questions left open. But perhaps these open
questions are of a more degree-theoretic flavor.  We will list three
questions here.

\begin{question}[Completeness]
  Which c.e.\ sets are automorphic to complete sets?
\end{question}

Of course, by \citet{Harrington.Soare:91}, we know that not every
c.e.\ set is automorphic to a complete set, and partial
classifications of precisely which sets can be found in
\citet{Downey.Stob:92} and Harrington and Soare
\cite{Harrington.Soare:96,Harrington.Soare:98}.

\begin{question}[Cone Avoidance]
  Given an incomplete c.e.\ degree $\mathbf{d}$ and an incomplete
  c.e.\ set $A$, is there an $\Ahat$ automorphic to $A$ such that
  $\mathbf{d} \not\leq_T \Ahat$?
\end{question}

\begin{question}[Can single jumps be coded into $\mathcal{E}$?]
  Let $J$ be C.E.A.\ in $\mathbf{0'}$ but not of degree
  $\mathbf{0''}$.  Is there a degree $\mathbf{a}$ such that
  $\mathbf{a'} \equiv_T J$ and, for all $A \in \mathbf{a}$, there is
  an $\Ahat$ with $A$ automorphic to $\Ahat$ and $\Ahat' <_T
  \mathbf{a'}$ or $\Ahat' |_T \mathbf{a'}$?
\end{question}

   \begin{question}[Can a single Turing degree be coded into
     $\mathcal{E}$?]
     Is there a degree $\mathbf{d}$ and an incomplete set $A$ such
     that, for all $\Ahat$ automorphic to $A$, $\mathbf{d} \leq
     \Ahat$?  $A \in \mathbf{d}$?
   \end{question}

   In a technical sense, these may not have a ``reasonable''
   answer. Thus the following seems a reasonable question.

\begin{question}
  Are these arithmetical questions?
\end{question}

In this paper we do not have the space to discuss the import of these
questions.  Furthermore, it not clear how this current work impacts
possible approaches to these questions.  At this point we will just
direct the reader to slides of a presentation of Cholak \cite{2006:c};
perhaps a paper reflecting on these issues will appear later.

One of the issues that will impact all of these questions are which
degrees can be realized in the orbits that we construct in
Theorem~\ref{main} and \ref{sec:maincor}.  A set is \emph{hemimaximal}
iff it is the nontrivial split of a maximal set.  A degree is
\emph{hemimaximal} iff it contains a hemimaximal set.
\citet{Downey.Stob:92} proved that the hemimaximal sets form an orbit,
and in some sense, this orbit is very large
degree-theoretically. While it is known by
\citet{Downey.Harrington:96} that there is no orbit containing sets of
all nonzero degrees, the orbit of hemimaximal sets contain
representatives of all jump classes (\citet{Downey.Stob:91}).

We are able to also show that we can construct our orbits to contain
at least a fixed hemimaximal degree (possibly along others) or contain
all hemimaximal degrees (again possibly along others).  However, what
is open is if every such orbit must contain a representative of every
hemimaximal degree or only hemimaximal degrees. For the proofs of
these claims, we direct the reader to Section~\ref{hemimaximal}.

\section{Past Work and Other Connections}

The paper \cite{cholak.downey.ea:nd} is a fourth paper in a series of
loosely connected papers, the previous three being by Harrington and
Cholak \cite{Cholak.Harrington:02}, \cite{Cholak.Harrington:03}, and
\cite{Cholak.Harrington:nd}.  We have seen above that results from
\cite{Cholak.Harrington:nd} determine the direction one must take to
prove Theorem~\ref{main}.  The above results from
\cite{Cholak.Harrington:nd} depend heavily on the main result in
\cite{Cholak.Harrington:03} whose proof depends on special
$\mathcal{L}$-patterns and several theorems about them which can be
found in \cite{Cholak.Harrington:02}. It is not necessary to
understand any of the above-mentioned theorems from any of these
papers to understand the proof of Theorem~\ref{main}.

But the proof of Theorem~\ref{main} does depend on Theorems~2.16,
2.17, and 5.10 of \citet{Cholak.Harrington:nd}; see
Section~\ref{sec:sketch-proof-theorem}. The proof of
Theorem~\ref{sec:maincor} also needs Theorem~6.3 of
\citet{Cholak.Harrington:nd}. The first two theorems are
straightforward but the third and fourth require work. The third is
another modified ``Extension Theorem.''  The fourth is what we might
call a ``Restriction Theorem''; it restricts the possibilities for
automorphisms.

Fortunately, we are able to use these four theorems from
\citet{Cholak.Harrington:nd} as {\em black boxes}.  These four
theorems provide a clean interface between the two papers.  If one
wants to understand the {\em proofs} of these four theorems one must
go to \citet{Cholak.Harrington:nd}; otherwise, the paper
\cite{cholak.downey.ea:nd} is completely independent from its three
predecessors.  In the next section we will explore the statements of
Theorems~5.10 and 6.3 of \citet{Cholak.Harrington:nd} in more detail.

\subsection{An Algebraic Framework}
\label{horse}

Crucial to the Theorem~\ref{main} is the following theorem of the
second two authors which demonstrates that $\Delta_3^0$ Extension
Lemmas are central to our understanding of the automorphism group of
${\mathcal E}^*$.

\begin{theorem}[Theorem~6.3 of \citet{Cholak.Harrington:nd}]\label{old2}
  Assume $D$ and $\Dhat$ are automorphic via $\Psi$. Then $D$ and
  $\Dhat$ are automorphic via $\Theta$ where $\Theta \upharpoonright
  \E(D)$ is $\Delta^0_3$. 
\end{theorem}

Theorem \ref{old2} says that {\em inside} any automorphism can be
thought of as $\Delta_3^0$.  The proof of this result and others we
will need relies heavily on the framework of the second two authors
who have recast the idea of an Extension Lemma {\em algebraically} so
that the dynamic notions of entry states and matching are replaced by
{\em extendible Boolean algebras} and {\em supports}.  In particular,
these proofs relies on Theorem \ref{theext}.  We will briefly discuss
these methods of Cholak and Harrington.

Fix a c.e.\ set $A$. Then the structure $\S(A)=\{B:\exists C(C\sqcup
B=A)\},$ the Boolean algebra of (c.e.)  {\em splits} of $A$.  Let
$\R(A)=\{R: R \subseteq A$ and $R$ computable$\},$ with $\S_R(A)$ the
quotient of $\S(A)$ by $\R(A),$ and $=^R,\subseteq^R$ the
corresponding quotient relations.  It is proven in Cholak and
Harrington \cite{Cholak.Harrington:nd} that $\S_R(A)$ is always a
$\Sigma_3^0$ Boolean algebra.  Ones that have representations of low
complexity are especially important to us.  A uniformly computable
listing $\S=\{S_i:i\in \omega\}$ of splits of $A$ is called an {\em
  effective listing of splits} of $A$ iff there exists another
uniformly computable listing of splits of $A$, $\{\hat{S}_i:i\in
\omega\}$, with $S_i\sqcup \hat{S}_i=A$ for all $i$.  The idea here is
based around the fact that those elements $x$ in $W_e$ and then later
enter $A$ give rise to effective listings of splits.

\begin{definition}[Cholak and Harrington \cite{Cholak.Harrington:nd}]
  A $\Sigma_3^0$ subalgebra $\B$ of $\S_R(A)$ is called {\em
    extendible} iff there exists a representation $\S$ and $B$ of $\B$
  such that $\S$ is an effective listing of splits of $A$ and $B$ is a
  $\Delta_3^0$ set.
\end{definition}

Again following Cholak and Harrington \cite{Cholak.Harrington:nd}, we
consider a partial map $\Theta$ between splits of $A$ and splits of
$\Ahat$ (for general sets $A$ and $\Ahat$) to be an isomorphism
between a substructure $\B$ of $\S_R(A)$ and a substructure $\Bhat$ of
$\S_R(\Ahat)$, iff $\Theta$ preserves $\subseteq^R$, for each
equivalence class $S_R$ of $\B$, if $S\in S_R$ then $\Theta(S)$
exists, and for each equivalence class $\hat{\S}_R$ of $\Bhat$,
$\Theta^{-1}(\Shat)$ exists for all $\Shat \in \hat{\S}_R$. Then two
extendible algebras $\B$ and $\Bhat$ are {\em extendibly isomorphic}
via $\Theta$ iff
\begin{enumerate}
\item There is an effective listing $\S$ and $B$ witnessing that $\B$
  is an extendible algebra.
\item There is an effective listing $\Shat$ and $\Bhat$ witnessing
  that $\hat{\B}$ is an extendible alga.
\item For all $i\in B$, there is a $j\in \Bhat$, with
  $\Theta(S_i)=\Shat_j,$ and
\item For each $j\in \Bhat,$ there is an $i\in B$ with
  $\Theta(S_i)=\Shat_j.$
\item The partial map $\Theta'$ indiced by $\Theta$ describes an
  isomorphism between $\B$ and $\hat{\B}$, as above.
\end{enumerate}

The first algebraic version of the Extension Lemma is the following.

\begin{theorem}[Cholak and Harrington \citet{Cholak.Harrington:nd}]
  \label{12}
  Let $\B\subseteq \S_R(A)$ and $\hat{\B}_R\subseteq \S_R(\Ahat)$ be
  two extendible Boolean algebras, which are $\Delta_3^0$ extendibly
  isomorphic via $\Theta$. Then there is a $\Phi$ which is a
  $\Delta_3^0$ isomorphism between $\E^*(A)$ and $\E^*(\Ahat)$, such
  that for all $i\in B$, $\Phi(S_i)=_R \Theta(S_i)$, and for all $i\in
  \Bhat,$ $\Phi^{_1}(\Shat_i)=_R \Theta^{-1}(\Shat_i).$
\end{theorem}

The key idea here is that it is possible to extend the extendible
isomorphism between $\B$ and $\hat{\B}$ to an isomorphism between
$\E^*(A)$ and $\E^*(\Ahat)$.  Of course, nothing comes without price,
and the {\em proof} of this (and similar) results, rely on dynamic
extension lemmas, of one type or another. In the paper Cholak and
Harrington \cite{Cholak.Harrington:nd}, Theorem \ref{12} is proven
using a modification of Cholak's {\em Translation Theorem}, Cholak
\cite{Cholak:94*1}

Whilst it is not directly pertinent to the present paper, we point out
how Cholak and Harrington applied theorems like Theorem \ref{12} using
the idea of supports.  This notion is related to the relationship
between $\L^*(A)$ and $\B$, in some sense focusing on the relationship
between the outside and the inside.  We say a c.e.\ set $S$ {\em
  supports} $X$ if $S\subseteq X$ and $(X-A)\sqcup S$ is c.e..  For
example, $W_e\searrow A$, the elements of $W_e$ which begin outside of
$A_s$ and then later enter $A$ support $W_e$. More generally, an
extendible algebra $\B$ supports a substructure $\L$ of $\L^*(A)$ (a
subcollection of $\{W_e\cup A,\cup\}$ modulo finite sets), if for all
$W\in \L$ there is an $i\in B$ with $S_i$ supporting $W$.

\begin{definition}[Cholak and Harrington \citet{Cholak.Harrington:nd}]
  Assume that $\L^*(A)$ and $\L^*(\Ahat)$ are isomorphic via $\Psi$,
  $\B$ and $\hat{\B}$ are extendible algebras isomorphic via $\Theta$,
  $\B$ supports $\L$ and $\hat{\B}$ supports $\hat{\L}$.  Then we say
  that $\Psi$ and $\Theta$ {\em preserve the supports } if for $W\in
  \L$, there is an $i\in B$ such that $S_i $ supports $W$ and
  $(\Psi(W\cup A)-\Ahat)\sqcup \Theta(S_i)$ is c.e., and for all
  $\What\in \hat{\L},$ there is an $i\in \hat{\B}$ with $\hat{S}_i$
  supporting $\What$ and $\Psi^{-1}(\What\cup\Ahat)-A) \sqcup
  \Theta^{-1}(\hat{S}_i) $ is c.e..\end{definition}

Then one algebraic version of the Extension Lemma is the following.

\begin{theorem}[Cholak and Harrington \cite{Cholak.Harrington:nd}]
  \label{theext}
  Assume that $\L^*(A)$ and $\L^*(\Ahat)$ are isomorphic via $\Psi$,
  $\B$ and $\hat{\B}$ are extendible algebras isomorphic via $\Theta$,
  $\B$ supports $\L^*(A)$ and $\hat{\B}$ supports $\L^*(\Ahat)$, with
  $\Psi$ and $\Theta$ preserving supports.

  Then there is an automorphism $\Lambda$ of $\E^*$ with
  $\Lambda(A)=\Ahat$, $\Lambda\upharpoonright \L^*(A)=\Psi$, and such
  that $\Lambda\upharpoonright \E^*(A)$ is $\Delta_3^0$.\end{theorem}

For reasons which become clear later, one final result from Cholak and
Harrington's paper we will need concerns extendible algebras of
computable sets.  An extendible algebra $\B$ of $\S_R(\omega)$ is
called a {\em extendible algebra of computable sets}, as the splits of
$\omega$ are computable sets.

\begin{theorem}[Theorem~5.10 of
  \citet{Cholak.Harrington:nd}]\label{interface2}
  Let $\B$ be an extendible algebra of computable sets and similarly
  for $\widehat{\B}$.  Assume the two are extendibly isomorphic via
  $\Pi$.  Then there is a $\Phi$ such that $\Phi$ is a $\Delta^0_3$
  isomorphism between $\E^*(A)$ and $\E^*(\widehat{A})$, $\Phi$ maps
  computable subsets to computable subsets, and, for all $R \in \B$,
  $(\Pi(R) - \widehat{A}) \sqcup \Phi(R \cap A)$ is computable\ (and
  dually).
\end{theorem}


\subsection{Some Algebraic Orbits}\label{sec:an-algebraic-orbit}

In \cite{Cholak.Harrington:nd}, Cholak and Harrington use Theorem
\ref{theext} to give algebraic proofs of many known theorems from the
literature such as the maximal and hemimaximal results. (Additionally
they use the algebraic methods for new results such as new orbits.)
We remark that the algebraic view does go back to Herrmann's proof
that a certain class of sets (now called {\em Hermann sets}) were
automorphic.  This result was proven using a hitherto unobserved {\em
  algebraic} consequence of the original Soare paper about the
preservation of computable sets under Soare's construction. Given that
is not well understood and is a critical proof of the
Theorem~\ref{main}, we would like to explore these algebraic proofs
with some more detail.

\begin{definition}
  $\C(A)$ is the set of $W_e$ such that either $\overline{A} \subseteq
  W_e$ or $ W_e \subseteq^* A$.
\end{definition}


\begin{theorem}[Soare's Automorphism Theorem \cite{Soare:74}]
  \label{apcomputable}
  Let $A$ and $\Ahat$ be two noncomputable \ce sets.
  \begin{enumerate}
  \item Then there is a $\Delta^0_3$ isomorphism $\Lambda$ between
    $\E(A) \cup \C(A)$ and $\E(\Ahat) \cup \C(\Ahat)$. Furthermore a
    $\Delta^0_3$-index for $\Lambda$ can be found uniformly from
    indexes for $A$ and $\Ahat$.
  \item In addition, $\Lambda$ preserves the computable subsets of
    $A$.
  \end{enumerate}
\end{theorem}

\citet{Soare:74} explicitly stated Theorem~\ref{apcomputable}.1.
Theorem~\ref{apcomputable}.2 was observed, in unpublished work, by
Herrmann.  Assume that $R$ is a computable subset of $A$. Herrmann's
observation was that $\overline{R} \in \C(A)$ and hence $\Lambda(R)
\sqcup \Lambda(\overline{R}) =^* \widehat{\omega}$ and therefore
$\Lambda$ maps $R$ to a computable subset of $\Ahat$.  This
observation of Herrmann was never published and is one of the key
facts he used in showing that the Herrmann sets form an orbit; see
\citet{Cholak.Downey.ea:01}.

\begin{theorem}[\citet{Soare:74}]
  The maximal sets form an orbit.
\end{theorem}

\begin{proof}
  Assume that $A$ and $\Ahat$ are maximal.  Then $\C(A) = \E$.  If $W
  \subseteq A$ then let $\Psi(W) = \Lambda(W)$. If $W \cup A =^*
  \omega$ there is a computable set $R_W$ such that $R_w \subseteq^*
  A$ and $\overline{R_W} \subseteq^* A$ and then let $\Psi(W) =
  \Lambda(W \cup R_W) \sqcup \overline{\Lambda(R_W)}$.  It is not
  difficult to show $\Psi$ is an automorphism.
\end{proof}

Recall that set is \emph{hemimaximal} iff it is the nontrivial split
of a maximal set.

\begin{theorem}[\citet{Downey.Stob:92}]
  The hemimaximal sets form an orbit.
\end{theorem}

\begin{proof}
  Assume $A_1\sqcup A_2 = A$ where the $A_i$s are not computable and
  $A$ is maximal.  Dually for $\Ahat$. Assume that $\Theta_i$ is an
  isomorphism from $\E^*(A_i)$ to $\E^*(\Ahat_i)$ that preserves the
  computable subsets (from Theorem~\ref{apcomputable}).

  As with the maximal sets, it is enough to define an isomorphism
  $\Lambda$ between $\E^*(A)$ and $\E^*(\Ahat)$ preserving the
  computable subsets. If $X \subseteq^* A$ then let $\Lambda(X) =
  \Theta_1(X \cap A_1) \sqcup \Theta_2(X \cap A_2)$.  Let ${R} \in
  {\R}(A)$. Then ${R} \cap A_i$ is computable.  So $\Theta_i({R} \cap
  A_i)$ is computable.  Hence $\Theta_1({R} \cap A_1) \sqcup
  \Theta_2({R} \cap A_2)$ is computable.  The complexity of the
  resulting automorphism is $\Delta^0_3$.
\end{proof}

\begin{definition}
  We say that a c.e.\ set $H$ is {\em strongly $r$-separable} if, for
  all c.e.\ sets $W$ disjoint from $H$, there is a computable set $R$
  such that $W\subset R$, $H\subset \overline{R}$, and $R-W$ is
  infinite.  We say that a set $H$ is {\em Herrmann} if it is both
  $\mathcal{D}$-maximal and strongly $r$-separable.
\end{definition}

\begin{theorem}[Herrmann, see
  \citet{Cholak.Downey.ea:01}]\label{sec:an-algebraic-orbit-1}
  The Herrmann sets form an orbit (under $\Delta^0_3$ automorphisms).
\end{theorem}

\begin{proof}
  Let $H$ be a Herrmann set.  Since $H$ is $\mathcal{D}$-maximal for
  all $W$ there is a $W_W$ such that either $W \subseteq^* H \sqcup
  W_W$ or $\overline{W} \subseteq^* H \sqcup W_W$.  Furthermore, since
  $H$ is Herrmann, for all $W$, there is a computable $R_W$ such that
  either $W \subseteq^* H \sqcup R_W$ or $\overline{W} \subseteq^* H
  \sqcup R_W$.  Note that finding $R_W$ and determining which case
  holds can be done using an oracle computable in $\mathbf{0}''$.

  Assume $\overline{W} \subseteq^* H \sqcup R_W$. Then $W \cup H
  \sqcup R_W =^* \omega$.  Therefore $(W \cap \overline{R_W}) \cup H
  =^* \overline{R}_W$. Recall that $X \backslash Y =\{ x | \exists s
  (x \in X_s-Y_s)\}$.  Hence $((W \cap \overline{R}_W)\backslash H)
  \sqcup (H \backslash (W \cap \overline{R}_W)) =^* \overline{R}_{W}$.
  Thus there is a computable subset $R_{H,W} = H \backslash (W \cap
  \overline{R}_W)$ of $H$ such that
  \begin{equation}\label{eq:22}
    W=^* (W \cap R_{H,W}) \sqcup (\overline{R}_{H,W}
    \cap \overline{R}_W)  \sqcup (W \cap R_W). 
  \end{equation}
  Again note that find $R_{H,W}$ can be using an oracle computable in
  $\mathbf{0}''$.

  Now using $\mathbf{0}''$ find a pairwise disjoint collection of
  $R_i$ such that $R_i \cap H = \emptyset$ and, for all $e$, $R_{W_e}
  \subseteq \bigsqcup_{i \leq e} R_e$. Since $H$ is Herrmann it is
  possible to find such a collection.  Do the same for $\widehat{H}$.

  Use Theorem~\ref{apcomputable} get $\Lambda$ mapping $\E(H)$ to
  $\E(\Hhat)$.  Let $p_i$ be a computable one-to-one onto map from
  $R_i$ to $\Rhat_i$ and, for $W \subseteq^*R_i$, let $\Lambda_i(W) =
  p_i(W)$. All of these maps take computable subsets to computable
  subsets.

  Now we will work on defining our automorphism $\Phi$. First assume
  $W \subseteq^* H \sqcup \bigsqcup R_i$. Let $\Phi(W) = \Lambda(H
  \cap W) \sqcup \bigsqcup \Lambda_i(R_i \cap W)$. It is not hard to
  show for $\Phi$, as defined so far, is order-preserving, for all
  $W$, $\Phi(W)$ is an r.e.\ set, and if $W$ is computable so is
  $\Phi(W)$.

  Now consider the case when $W \not\subseteq^* H \sqcup \bigsqcup
  R_i$.  In that case, Equation~\ref{eq:22} holds, and we can use that
  to define $\Phi(W)$ in terms of subsets of $H \sqcup \bigsqcup R_i$:
  \begin{equation*}\label{eq:2}
    \Phi(W)=^* \Phi(W \cap R_{H,W}) \sqcup (\overline{\Phi(R_{H,W})}
    \cap \overline{\Phi(R_W)})  \sqcup \Phi(W \cap R_W). 
  \end{equation*}
  As defined $\Phi(W)$ is an r.e.\ set.  It is not difficult show
  $\Phi$ is order preserving and hence well-defined. Thus $\Phi$ is an
  automorphism taking $H$ to $\Hhat$.
\end{proof}

For more on preserving the computable sets and an algebraic proof of
Theorem~\ref{apcomputable} we refer the reader to
\citet{Cholak.Harrington:nd}.

\section{A Sketch of the Proof of Theorem~\ref{main}}

\label{sec:sketch-proof-theorem}
The proof of Theorem \ref{main} is quite complex and involves several
ingredients. The proof will be easiest to understand if we introduce
each of the relevant ingredients in context.

The following theorem will prove be to useful.

\begin{theorem}[Folklore\footnote{See Section~5 of the full paper for
    more information and a proof.}]\label{folkloreI}
  There is a computable listing $T_i$ of computable infinite branching
  trees and a computable infinite branching tree $T_{\Sigma^1_1}$ such
  that the set $\{ i : T_{\Sigma^1_1} \cong T_i\}$ is
  $\Sigma^1_1$-complete.
\end{theorem}

The idea for the proof of Theorem~\ref{main} is to code each of the
above $T_i$s into the orbit of $A_{T_i}$. Informally let
$\mathcal{T}(A_T)$ denote this encoding; $\mathcal{T}(A_T)$ will not
be defined in this announcement. But we will discuss in some details
some of the ingredients and resulting complexity.  The game plan is as
follows:
  
\begin{enumerate}
\item \textbf{Coding:} For each $T$ build an $A_T$ such that $T \cong
  \mathcal{T}(A_T)$ via an isomorphism $\Lambda\leq_T
  \bf{0}^{(2)}$. 
\item \textbf{Coding is preserved under automorphic images:} If $\Ahat
  \approx A_T$ via an automorphism $\Phi$ then $\mathcal{T}(\Ahat)$
  exists and $\mathcal{T}(\Ahat) \cong T$ via an isomorphism
  $\Lambda_{\Phi}$, where $\Lambda_\Phi \leq_T \Phi \join
  \bf{0}^{(2)}$.
\item \textbf{Sets coding isomorphic trees belong to the same orbit:}
  If $T \cong \widehat{T}$ via isomorphism $\Lambda$ then $A_T \cong
  A_{\widehat{T}}$ via an automorphism $\Phi_\Lambda$ where
  $\Phi_\Lambda\leq_T \Lambda \join \bf{0}^{(2)}$.
\end{enumerate}

Thus $A_{T_{\Sigma^1_1}}$ and $A_{T_i}$ are in the same orbit iff
$T_{\Sigma^1_1}$ and $T_i$ are isomorphic.  Since the latter question
is $\Sigma^1_1$-complete so is the former question.

We will build a pairwise disjoint collection of $D_\chi$s to code the
tree $T$.  $A_T = D_\lambda$ will code the empty node in $T$.  The
basic module to construct an $D$ involves the construction of a
computable set $R$ and a subset $M$ of $R$. There will be infinitely
many pairwise disjoint $R$s. Inside $R$s the $D$s will be Friedberg
splits of $M$.  It is well understood how to split an r.e.\ set $M$
into Friedberg splits.

Depending on the construction either $M=^*R$ or $M$ is maximal inside
$R$; i.e.\ $M \cup \overline{R}$ is maximal.  If $M$ is maximal inside
$R$ then $D$ is hemimaximal inside $R$ and we say $D$ \emph{lives} in
$R$. If $M =^* R$ then, inside $R$, $D$ is a computable set and $D$
does not live in $R$. Note that $D$ living in $R$ is a definable
property as is $D$ is computable in $R$.

The big issue of the construction will be to decide when $D$ lives in
$R$ and when not.  It is well understood how to construct an maximal
set inside $R$.  Furthermore it also well understood how use a dumping
argument to alter the maximal set construction to force $M=^*R$. We
will use a dumping construction to alter the construction of the
desired $M$.  The decision whether $D$ lives in $R$ will be handled by
a tree argument.  Here we will not discuss the tree argument but
discuss issues that go into deciding whether $D$ lives in $R$.

Each node $\chi$ of $T$ will have infinitely many pairwise disjoint
$R_{\chi,i}$ associated with $\chi$.  For these $R_{\chi,i}$ the
corresponding set is $M_{\chi,i}$.  $D_\chi$ will always be a
Friedberg split of $M_{\chi,i}$.  If $\chi^+$ is a successor of $\chi$
in $T$ then, for almost all $R_{\chi,i}$, $D_{\chi^+}$ be a Friedberg
split of $M_{\chi,i}$.  The collection of the all $R_{\chi,i}$ will be
pairwise disjoint.  When constructed in this fashion the $D_\chi$ code
$T$ and hence we have part 1 of the game plan under control.



We have to work on part 2 and 3 of the game plan. For part 3 we are
going to set things up so that if we know where the $D_{\chi}$ and
$R_{\chi,i}$ go then we will be able to construct the desired
automorphism.

Lets look at the hatted side of the construction briefly to work on
Part 2.  We are told $A_T$ goes to $\Ahat$.  We can gather together a
collection of pairwise disjoint computable sets $\Rhat_{\lambda,i}$
such that either $\Ahat$ lives in $\Rhat_{\lambda,i}$ or $\Ahat$ is
computable inside $\Rhat_{\lambda,i}$.  We can assume that this list
is maximal; i.e.\ if $\Ahat$ lives in $\Rhat$ or is computable in
$\Rhat$, then $\Rhat \subseteq \bigsqcup \Rhat_{\lambda,i}$.  There
may many such lists. But we can show modulo a computable set each
$R_{\lambda,i}$ must be sent to some $\Rhat_{\lambda,j}$.

Now we are in position to pick out the successors of $\Ahat$ on the
hatted side.  They are the sets $\Dhat_{\lambda^+}$ which live into
almost all of $\Rhat_{\lambda,i}$ that $\Ahat$ lives in.  More or less
these sets must be the automorphic images of the
$D_{\lambda^+,i}$s. We better ensure that these sets behave like a
good successor.  The only way to do is to control how the
$D_{\lambda^+,i}$s behave.

Any set $D$ which looks like a $D_\chi$ must be split of one of the
$D_\chi$s.  If, for infinitely many $\chi$ and $i$, $D$ lives in
$R_{\chi,i}$ then there must be exactly one $\chi$ such that, for all
$i$, $D$ lives in $R_{\chi,i}$ iff $D_\chi$ lives in $R_{\chi,i}$ and
for almost all $i$, $D$ lives in $R_{\chi^-,i}$ iff $D_{\chi^-}$ lives
in $R_{\chi^-,i}$.  In this case $D$ will be a Friedberg split of
$D_\chi$.

If we can get the $D$ to behave properly (as discussed above) than
their automorphic images must also behave properly. So the
$\Dhat_{\lambda^+}$ which we have found above in fact code the
successors of $\lambda$ in $\That$.  Once we have this we can find
their successors.

For each $\lambda^+ = \widehat{\chi}$ we can find a collection of
pairwise disjoint computable sets $\Rhat_{\widehat{\chi},i}$ such that
either $\Dhat_{ \widehat{\chi}}$ lives in $\Rhat_{\widehat{\chi},i}$
or $\Dhat_{\widehat{\chi}}$ is computable inside
$\Rhat_{\widehat{\chi},i}$.  We can assume that this list is maximal;
i.e.\ if $\Dhat_{ \widehat{\chi}}$ lives in $\Rhat$ or is computable
in $\Rhat$, then $\Rhat \subseteq \bigsqcup \Rhat_{\widehat{\chi},i}
\sqcup \bigsqcup \Rhat_{\lambda,i}$.  There may many such lists. Given
that $\Dhat_{ \widehat{\chi}}$ does code a successor of $\Ahat$, there
is some node $\chi$ of length $1$ such that modulo a computable set
each $R_{\chi,i}$ must be sent to some $\Rhat_{\widehat{\chi},j}$.
Now we can bootstrap our way to find successor of $\Dhat_{
  \widehat{\chi}}$ and so on.  Hence part 2 is now under control.

So, for each possible $D$, we will make sure that if either $D$ is
hemimaximal inside finite many $R_{\chi,i}$ or $D$ behaviors as above;
i.e.\ $D$ will be a Friedberg split of some $D_\chi$. To do this we
will use a coherence/state argument not unlike the argument used to
constructed maximal sets.  The state of $R_{\chi,i}$ will be those $D$
which $D$ is hemimaximal in $R_{\chi,i}$. Determining the state of a
set is $\Sigma^0_3$ rather than $\Sigma^0_1$.  Hence this is another
reason we must do this whole construction on a tree.  If $R_{\chi,i}$
is in an low $e$-state then we must dump it. We dump $R_{\chi,i}$ by
making $M_{\chi,i} =^* R_{\chi,i}$ as discussed above.

As for the collection of all $R_{\chi,i}$ we want them to have the
property that for all $W$, either $W$ is a split of the $D_\chi$s or
there is a finite set $F$ of $\chi$ and $i$ such that either $W
\subseteq^* \bigsqcup_{(\chi,i) \in F} R_{\chi,i}$ or $\overline{W}
\subseteq^* \bigsqcup_{(\chi,i) \in F} R_{\chi,i}$. This cannot be
achieved via an effective construction but is achievable on a tree.
This has the side effect that there will be more sets $R_{\chi,i}$
where $M_{\chi,i} =^* R_{\chi,i}$.

We will also construct that $R_{\chi,i}$ such that for all $\chi$, the
$R_{\chi^-,i}$ and $R_{\chi,i}$ form an extendible algebra of
computable sets, $\mathcal{B}_\chi$.  We wish to use
Theorem~\ref{interface2}. For each $\chi$, the map $R_{\chi,i}$ to
$\Rhat_{\widehat{\chi},i}$ is an extendible isomorphism. We want to
claim that we can use these pieces and an isomorphism between $T$ and
$\That$ to construct an automorphism not unlike what we did in
Section~\ref{sec:an-algebraic-orbit} and, in particular, with the
proof of Theorem~\ref{sec:an-algebraic-orbit-1}.  With one caveat this
is the case.

The caveat is that if $\chi$ and $\widehat{\chi}$ have the same length
we want to ensure that $D_\chi$ lives in $R_{\chi,i}$ iff
$\Dhat_{\widehat{\chi}}$ lives in $\Rhat_{\widehat{\chi},i}$.  This
requires two more additions to the above description.  First we must
construct all $A_T$s using the same tree construction.  Second the
coherence/state argument above must be extended so that we the above
homogeneous is preserved.  So if $M_{\chi,i} =^* R_{\chi,i}$ then for
all $\widehat{\chi}$ of the same length as $\chi$,
$M_{\widehat{\chi},i} =^* R_{\widehat{\chi},i}$.

\subsection{Invariants}

\label{sec:last}

It might appear that $\mathcal{T}(A)$ is an invariant which determines
the orbit of $A$. But there is no reason to believe for an arbitrary
$A$ that $\mathcal{T}(A)$ is well defined.  The following theorem
shows that $\mathcal{T}(\Ahat)$ is an invariant as far as the orbits
of the $A_T$s are concerned.  In the full paper, we prove a more
technical version of the following theorem.

\begin{theorem}\label{sec:invariant}
  If $\Ahat$ and $A_T$ are automorphic via $\Psi$ and $T \cong
  \mathcal{T}(\Ahat)$ via $\Lambda$ then $A_T \approx \Ahat$ via
  $\Phi_\Lambda$ where $\Phi_\Lambda\leq_T \Lambda \join
  \bf{0}^{(8)}$.
\end{theorem}

\begin{proof}[Sketch]
  For $A_T$ the above construction gives us a $\mathbf{0''}$ (they are
  constructed on the true path) listing of the sets $D_\chi$,
  $R_{\chi,i}$, and $M_{\chi,i}$. So they are available for us to use
  here. The idea is to recover images of these sets on the hatted
  side.  This recovery relies on Theorem~6.3 of
  \citet{Cholak.Harrington:nd} that if $D$ and $\Dhat$ are automorphic
  via $\Psi$, then $D$ and $\Dhat$ are automorphic via $\Theta$ where
  $\Theta \upharpoonright \E(D)$ is $\Delta^0_3$, as well as more
  intricate material on extendible algebras (in particular, careful
  application of Theorem 5.10 of \citet{Cholak.Harrington:nd}).  Then
  using these recovered sets we constructed the desired automorphism
  as hinted to above.  It takes $\Lambda \join \bf{0}^{(8)}$ to
  recover the needed sets.  The construction of the automorphism needs
  the recovered sets and an oracle for $\bf{0}^{(2)}$.
\end{proof}

\subsection{Properly $\Delta^0_\alpha$ orbits}

\begin{theorem}[Folklore\footnote{See Section~5 of the full paper
    \cite{cholak.downey.ea:nd} for more information and a
    proof.}]\label{sec:folkloredark}
  For all 
  finite $\alpha$ there is a computable tree $T_{i_\alpha}$ from the
  list in Theorem~\ref{folkloreBA} such that, for all computable trees
  $T$, $T$ and $T_{i_\alpha}$ are isomorphic iff $T$ and
  $T_{i_\alpha}$ are isomorphic via an isomorphism computable in
  $\text{deg}(T) \join 0^{(\alpha)}$.  But, for all $\beta < \alpha$
  there is an $i^*_\beta$ such that $T_{i^*_\beta}$ and $T_{i_\alpha}$
  are isomorphic but are not isomorphic via an isomorphism computable
  in $0^{(\beta)}$.
\end{theorem}

It is open if the above theorem holds for all $\alpha$ such that
$\omega \geq \alpha < \wock$.  But if it does then so does the theorem
below.

\begin{theorem}
  For all finite $\alpha > 8$ there is a properly $\Delta^0_{\alpha}$
  orbit.
\end{theorem}

\begin{proof}
  Assume that $A_{T_{i_\alpha}}$ and $\Ahat$ are automorphic via an
  automorphism $\Phi$.  Hence, by part 2 of the game plan,
  $\mathcal{T}(\Ahat)$ and $T_{i_\alpha}$ are isomorphic.  Since
  $\mathcal{T}(\Ahat)$ is computable in $0^{(8)}$, $\alpha > 8$, and
  by Theorem~\ref{sec:folkloredark}, $\mathcal{T}(\Ahat)$ and
  $T_{i_\alpha}$ via a $\Lambda \leq_T 0^{(\alpha)}$.  By
  Theorem~\ref{sec:invariant}, $\Ahat$ and $A_{T_{i_\alpha}}$ are
  automorphic via an automorphism computable in $0^{(\alpha)}$.

  Fix $\beta$ such that $8 \geq \beta < \alpha$.  By part 3 of the
  game plan and the above paragraph, $A_{T_{i_\alpha}}$ and
  $A_{T_{i^*_\beta}}$ are automorphic via an automorphism computable
  in $0^{(\alpha)}$.  Now assume $A_{T_{i^*_\beta}} \approx
  A_{T_{i_\alpha}}$ via $\Phi$.  By part 2 of the game plan,
  $\mathcal{T}(A_{T_{i^*_\beta}}) \cong T_{i_\alpha}$ via
  $\Lambda_{\Phi}$, where $\Lambda_\Phi \leq_T \Phi \join
  \bf{0}^{(2)}$.  Since $\mathcal{T}(A_{T_{i^*_\beta}})$ is computable
  in $0^{(8)}$ and $\mathcal{T}(A_{T_{i^*_\beta}})$ is isomorphic to
  $T_{i^*_\beta}$ via an isomorphism computable in $0^{(\beta)}$ (part
  1 of the game plan), by Theorem~\ref{sec:folkloredark},
  $\Lambda_\Phi >_T 0^{(\beta)}$.  Hence $\Phi >_T 0^{(\beta)}$.
\end{proof}

\subsection{Our Orbits and Hemimaximal Degrees}\label{hemimaximal}

Recall that set is \emph{hemimaximal} iff it is the nontrivial split
of a maximal set.  A degree is \emph{hemimaximal} iff it contains a
hemimaximal set.

Let $T$ be given.  Construction $A_T$ as above.  For all $i$, either
$A_T$ is hemimaximal in $R_i$ or $A_T \cap R_i$ is computable. If
$A_T$ is hemimaximal in $R_i$ then $A_T \cap R_i$ is a split of
maximal set $M \sqcup \overline{ R}_i$ and hence $A_T = (A_T \cap
R_i)$ is a hemimaximal set.  $A_T = \bigsqcup_{i \in \omega} (A_T \cap
R_i)$ where $A_T \cap R_i$ is either hemimaximal or computable.  So
the degree of $A_T$ is the infinite join of hemimaximal degrees.  It
is not known if the (infinite) join of hemimaximal degrees is
hemimaximal.  Moreover, this is not an effective infinite join. But if
we control the degrees of $A_T \cap R_i$ we can control the degree of
$A_T$.  By modifying our proofs we can achieve the following degree
controls.

\begin{theorem}\label{AThemi}
  Let $H$ be hemimaximal.  We can construct $A_T$ such that $A_T
  \equiv_T H$.
\end{theorem}

Indeed, we can show the following.

\begin{theorem}
  There is an $A_T$ whose orbits contain a representative of every
  hemimaximal degree, and hence of all jump classes.
\end{theorem}


\begin{thebibliography}{39}
\providecommand{\natexlab}[1]{#1}
\providecommand{\url}[1]{\texttt{#1}}
\expandafter\ifx\csname urlstyle\endcsname\relax
  \providecommand{\doi}[1]{doi: #1}\else
  \providecommand{\doi}{doi: \begingroup \urlstyle{rm}\Url}\fi

\bibitem[Ash and Knight(2000)]{MR1767842}
C.~J. Ash and J.~Knight.
\newblock \emph{Computable structures and the hyperarithmetical hierarchy},
  volume 144 of \emph{Studies in Logic and the Foundations of Mathematics}.
\newblock North-Holland Publishing Co., Amsterdam, 2000.
\newblock ISBN 0-444-50072-3.

\bibitem[Cholak(1994)]{Cholak:94*1}
P.~Cholak.
\newblock The translation theorem.
\newblock \emph{Arch. Math. Logic}, 33:\penalty0 87--108, 1994.

\bibitem[Cholak et~al.(1992)Cholak, Downey, and Stob]{Cholak.Downey.ea:92*1}
P.~Cholak, R.~Downey, and M.~Stob.
\newblock Automorphisms of the lattice of recursively enumerable sets: Promptly
  simple sets.
\newblock \emph{Trans. Amer. Math. Soc.}, 332:\penalty0 555--570, 1992.

\bibitem[Cholak(2006)]{2006:c}
Peter Cholak.
\newblock The {C}omputably {E}numerable {S}ets: the {P}ast, the {P}resent and
  the {F}uture.
\newblock Theory and Applications of Models of Computation, 2006, Beijing
  China, Slides can be found at \url{http://www.nd.edu/~cholak}, 2006.

\bibitem[Cholak(1991)]{Cholak:91}
Peter Cholak.
\newblock \emph{Automorphisms of the lattice of recursively enumerable sets}.
\newblock PhD thesis, University of Wisconsin, 1991.

\bibitem[Cholak(1995)]{Cholak:95}
Peter Cholak.
\newblock Automorphisms of the lattice of recursively enumerable sets.
\newblock \emph{Mem. Amer. Math. Soc.}, 113\penalty0 (541):\penalty0 viii+151,
  1995.
\newblock ISSN 0065-9266.

\bibitem[Cholak and Harrington(2003)]{Cholak.Harrington:03}
Peter Cholak and Leo Harrington.
\newblock Isomorphisms of splits of computably enumerable sets.
\newblock \emph{J. of Symbolic Logic}, 68\penalty0 (3):\penalty0 1044--1064,
  2003.

\bibitem[Cholak and Harrington()]{Cholak.Harrington:nd}
Peter Cholak and Leo~A. Harrington.
\newblock Extension theorems, orbits, and automorphisms of the computably
  enumerable sets.
\newblock To appear in Trans. Amer. Math. Soc. Final version as of 8/31/2005.
  math.LO/0408279.

\bibitem[Cholak et~al.()Cholak, Downey, and Harrington]{cholak.downey.ea:nd}
Peter Cholak, Rod Downey, and Leo~A. Harrington.
\newblock On the {O}rbits of {C}omputable {E}numerable {S}ets.
\newblock Submitted. math.LO/0607264.

\bibitem[Cholak et~al.(2001)Cholak, Downey, and Herrmann]{Cholak.Downey.ea:01}
Peter Cholak, Rod Downey, and Eberhard Herrmann.
\newblock Some orbits for {$\mathcal{E}$}.
\newblock \emph{Ann. Pure Appl. Logic}, 107\penalty0 (1-3):\penalty0 193--226,
  2001.
\newblock ISSN 0168-0072.

\bibitem[Cholak and Downey(2004)]{MR2005e:03092}
Peter~A. Cholak and Rod Downey.
\newblock Invariance and noninvariance in the lattice of {$\Pi\sp 0\sb 1$}
  classes.
\newblock \emph{J. London Math. Soc. (2)}, 70\penalty0 (3):\penalty0 735--749,
  2004.
\newblock ISSN 0024-6107.

\bibitem[Cholak and Harrington(2000)]{Cholak.Harrington:00}
Peter~A. Cholak and Leo~A. Harrington.
\newblock Definable encodings in the computably enumerable sets.
\newblock \emph{Bull. Symbolic Logic}, 6\penalty0 (2):\penalty0 185--196, 2000.
\newblock A copy can be found at \url{http://www.nd.edu/~cholak}.

\bibitem[Cholak and Harrington(2002)]{Cholak.Harrington:02}
Peter~A. Cholak and Leo~A. Harrington.
\newblock On the definability of the double jump in the computably enumerable
  sets.
\newblock \emph{J. Math. Log.}, 2\penalty0 (2):\penalty0 261--296, 2002.
\newblock ISSN 0219-0613.

\bibitem[Downey and Montalb{\'a}n(2006)]{Downey.Montalban:nd}
R.~Downey and A.~Montalb{\'a}n.
\newblock Slender classes.
\newblock Submitted, 2006.

\bibitem[Downey and Stob(1991)]{Downey.Stob:91}
R.~G. Downey and M.~Stob.
\newblock Jumps of hemimaximal sets.
\newblock \emph{Z. Math. Logik Grundlag. Math.}, 37:\penalty0 113--120, 1991.

\bibitem[Downey and Stob(1992)]{Downey.Stob:92}
R.~G. Downey and M.~Stob.
\newblock Automorphisms of the lattice of recursively enumerable sets: Orbits.
\newblock \emph{Adv. in Math.}, 92:\penalty0 237--265, 1992.

\bibitem[Downey and Harrington(1996)]{Downey.Harrington:96}
Rod Downey and Leo Harrington.
\newblock There is no fat orbit.
\newblock \emph{Ann. Pure Appl. Logic}, 80\penalty0 (3):\penalty0 277--289,
  1996.
\newblock ISSN 0168-0072.

\bibitem[Friedberg(1957)]{Friedberg:57}
R.~M. Friedberg.
\newblock A criterion for completeness of degrees of unsolvability.
\newblock \emph{J. Symbolic Logic}, 22:\penalty0 159--160, 1957.

\bibitem[Goncharov et~al.(2004)Goncharov, Harizanov, Knight, and
  Shore]{MR2058190}
Sergey~S. Goncharov, Valentina~S. Harizanov, Julia~F. Knight, and Richard~A.
  Shore.
\newblock {$\Pi^1_1$} relations and paths through {$\mathcal{O}$}.
\newblock \emph{J. Symbolic Logic}, 69\penalty0 (2):\penalty0 585--611, 2004.
\newblock ISSN 0022-4812.

\bibitem[Harrington and Soare(1998)]{Harrington.Soare:98}
Leo Harrington and Robert~I. Soare.
\newblock Codable sets and orbits of computably enumerable sets.
\newblock \emph{J. Symbolic Logic}, 63\penalty0 (1):\penalty0 1--28, 1998.
\newblock ISSN 0022-4812.

\bibitem[Harrington and Soare(1991)]{Harrington.Soare:91}
Leo~A. Harrington and Robert~I. Soare.
\newblock {P}ost's program and incomplete recursively enumerable sets.
\newblock \emph{Proc. Nat. Acad. Sci. U.S.A.}, 88:\penalty0 10242--10246, 1991.

\bibitem[Harrington and Soare(1996{\natexlab{a}})]{Harrington.Soare:96}
Leo~A. Harrington and Robert~I. Soare.
\newblock The ${\Delta}\sp 0\sb 3$-automorphism method and noninvariant classes
  of degrees.
\newblock \emph{J. Amer. Math. Soc.}, 9\penalty0 (3):\penalty0 617--666,
  1996{\natexlab{a}}.
\newblock ISSN 0894-0347.

\bibitem[Harrington and Soare(1996{\natexlab{b}})]{Harrington.Soare:96*1}
Leo~A. Harrington and Robert~I. Soare.
\newblock Definability, automorphisms, and dynamic properties of computably
  enumerable sets.
\newblock \emph{Bull. Symbolic Logic}, 2\penalty0 (2):\penalty0 199--213,
  1996{\natexlab{b}}.
\newblock ISSN 1079-8986.

\bibitem[Kleene and Post(1954)]{Kleene.Post:54}
Stephen~C. Kleene and Emil~L. Post.
\newblock The upper semi-lattice of degrees of recursive unsolvability.
\newblock \emph{Ann. of Math. (2)}, 59:\penalty0 379--407, 1954.

\bibitem[Lachlan(1968)]{Lachlan:68*3}
Alistair~H. Lachlan.
\newblock On the lattice of recursively enumerable sets.
\newblock \emph{Trans. Amer. Math. Soc.}, 130:\penalty0 1--37, 1968.

\bibitem[Lerman and Soare(1980)]{Lerman.Soare:80*1}
Manuel Lerman and Robert~I. Soare.
\newblock $d$-simple sets, small sets, and degree classes.
\newblock \emph{Pacific J. Math.}, 87\penalty0 (1):\penalty0 135--155, 1980.
\newblock ISSN 0030-8730.

\bibitem[Maass(1984)]{Maass:84}
W.~Maass.
\newblock On the orbit of hyperhypersimple sets.
\newblock \emph{J. Symbolic Logic}, 49:\penalty0 51--62, 1984.

\bibitem[Maass and Stob(1983)]{Maass.Stob:83}
W.~Maass and M.~Stob.
\newblock The intervals of the lattice of recursively enumerable sets
  determined by major subsets.
\newblock \emph{Ann. Pure Appl. Logic}, 24:\penalty0 189--212, 1983.

\bibitem[Marchenkov(1976)]{Marchenkov:76}
S.~S. Marchenkov.
\newblock A class of incomplete sets.
\newblock \emph{Math. Z.}, 20:\penalty0 473--487, 1976.

\bibitem[Martin(1966)]{Martin:66*1}
D.~A. Martin.
\newblock Classes of recursively enumerable sets and degrees of unsolvability.
\newblock \emph{Z. Math. Logik Grundlag. Math.}, 12:\penalty0 295--310, 1966.

\bibitem[Muchnik(1956)]{Muchnik:56}
A.~A. Muchnik.
\newblock On the unsolvability of the problem of reducibility in the theory of
  algorithms.
\newblock \emph{Dokl. Akad. Nauk SSSR}, N. S. 108:\penalty0 194--197, 1956.

\bibitem[Post(1944)]{Post:44}
Emil~L. Post.
\newblock Recursively enumerable sets of positive integers and their decision
  problems.
\newblock \emph{Bull. Amer. Math. Soc.}, 50:\penalty0 284--316, 1944.

\bibitem[Slaman and Woodin(1989)]{Slaman.Woodin:conjecture}
Theodore~A. Slaman and W.~Hugh Woodin.
\newblock Personal communication, 1989.

\bibitem[Soare(1974{\natexlab{a}})]{Soare:74}
Robert~I. Soare.
\newblock Automorphisms of the lattice of recursively enumerable sets {I}:
  maximal sets.
\newblock \emph{Ann. of Math. (2)}, 100:\penalty0 80--120, 1974{\natexlab{a}}.

\bibitem[Soare(1974{\natexlab{b}})]{Soare:74*1}
Robert~I. Soare.
\newblock Automorphisms of the lattice of recursively enumerable sets.
\newblock \emph{Bull. Amer. Math. Soc.}, 80:\penalty0 53--58,
  1974{\natexlab{b}}.

\bibitem[Turing(1939)]{Turing:39}
Alan~M. Turing.
\newblock Systems of logic based on ordinals.
\newblock \emph{Proc. London Math. Soc. (3)}, 45:\penalty0 161--228, 1939.

\bibitem[Weber(2004)]{Weber:04}
Rebecca Weber.
\newblock \emph{A definable relation between c.e.\ sets and ideals}.
\newblock PhD thesis, University of Notre Dame, 2004.

\bibitem[Weber(2006)]{Weber:06}
Rebecca Weber.
\newblock Invariance in {$\mathcal{E}^*$} and {$\mathcal{E}_{\Pi}$}.
\newblock \emph{Trans. Amer. Math. Soc.}, 358\penalty0 (7):\penalty0 3023--3059
  (electronic), 2006.
\newblock ISSN 0002-9947.

\bibitem[White(2000)]{walkerwhite:00}
Walker~M. White.
\newblock \emph{Characterizations for Computable Structures}.
\newblock PhD thesis, Cornell University, Ithaca, NY, USA, 2000.

\end{thebibliography}

\end{document}